\documentclass[a4]{amsart}

\input xypic
\usepackage{amssymb, amsmath} 

\newtheorem{theorem}{Theorem}[section]
\newtheorem{proposition}[theorem]{Proposition}
\newtheorem{lemma}[theorem]{Lemma}
\newtheorem{corollary}[theorem]{Corollary}

\theoremstyle{definition}

\newtheorem{remark}[theorem]{Remark}

\newcommand{\anick}{\ensuremath{T^{2n+1}(p^{r})}} 
\newcommand{\curly}{\ensuremath{S^{2n+1}\{p^{r}\}}} 
\newcommand{\jkm}{\ensuremath{J_{k-1}(S^{2n})}}
\newcommand{\jk}{\ensuremath{J_{k}(S^{2n})}} 
\newcommand{\pkm}{\ensuremath{P_{k-1}}}
\newcommand{\moore}{\ensuremath{P^{2n+1}(p^{r})}}

\newcommand{\hlgy}[1]{\ensuremath{H_{*}(#1)}} 
\newcommand{\rhlgy}[1]{\ensuremath{\widetilde{H}_{*}(#1)}} 

\newcounter{bean}

\newcommand{\namedright}[3]{\ensuremath{#1\stackrel{#2}
 {\longrightarrow}#3}}
\newcommand{\nameddright}[5]{\ensuremath{#1\stackrel{#2}
 {\longrightarrow}#3\stackrel{#4}{\longrightarrow}#5}} 
\newcommand{\namedddright}[7]{\ensuremath{#1\stackrel{#2}
 {\longrightarrow}#3\stackrel{#4}{\longrightarrow}#5
  \stackrel{#6}{\longrightarrow}#7}} 

\newcommand{\larrow}{\relbar\!\!\relbar\!\!\rightarrow} 
\newcommand{\llarrow}{\relbar\!\!\relbar\!\!\larrow}

\newcommand{\lnameddright}[5]{\ensuremath{#1\stackrel{#2}
 {\larrow}#3\stackrel{#4}{\larrow}#5}} 
 
\newcommand{\llnamedright}[3]{\ensuremath{#1\stackrel{#2}
 {\llarrow}#3}}
\newcommand{\llnameddright}[5]{\ensuremath{#1\stackrel{#2}
 {\llarrow}#3\stackrel{#4}{\llarrow}#5}}

\newcommand{\qqed}{\hfill\Box} 
 
\newcommand{\zmodp}{\ensuremath{\mathbb{Z}\mathit{/p}\mathbb{Z}}}

\begin{document} 

%%% Title 

\title{Anick's fibration and the odd primary homotopy 
       exponent of spheres}
\author{Stephen D. Theriault}
\address{Department of Mathematical Sciences, University
         of Aberdeen, Aberdeen AB24 3UE, United Kingdom}
\email{s.theriault@maths.abdn.ac.uk}


%%% Abstract 

\begin{abstract} 
For primes $p\geq 3$, Cohen, Moore, and Neisendorfer showed that 
the exponent of the $p$-torsion in the homotopy groups of $S^{2n+1}$ 
is $p^{n}$. This was obtained as a consequence of a thorough 
analysis of the homotopy theory of Moore spaces. Anick further developed  
this for $p\geq 5$ by constructing a homotopy fibration 
\(\nameddright{S^{2n-1}}{}{\anick}{}{\Omega S^{2n+1}}\) 
whose connecting map is degree~$p^{r}$ on the bottom cell. A much simpler  
construction of such a fibration for $p\geq 3$ was given by Gray and 
the author using new methods. In this paper the new methods are used 
to start over, first constructing Anick's fibration for $p\geq 3$, and 
then using it to obtain the exponent result for spheres. 
\end{abstract} 

\maketitle

\section{Introduction} 
\label{sec:intro} 

Cohen, Moore, and Neisendorfer's~\cite{CMN1,CMN2,N2} 
determination of the odd primary homotopy exponent of spheres 
is a miletone in homotopy theory. They showed that 
if $p$ is an odd prime then $p^{n}$ is the least power of $p$ 
which annihilates the $p$-torsion in $\pi_{\ast}(S^{2n+1})$. This 
was obtained as a consequence of a geometric result. 
Localize spaces and maps at $p$. Let 
\(E^{2}\colon\namedright{S^{2n-1}}{}{\Omega^{2} S^{2n+1}}\) 
be the double suspension. For each $n\geq 1$, they constructed a map 
\(\varphi\colon\namedright{\Omega^{2} S^{2n+1}}{}{S^{2n-1}}\)   
with the property that the composite 
\(\nameddright{\Omega^{2} S^{2n+1}}{\varphi}{S^{2n-1}}{E^{2}} 
   {\Omega^{2} S^{2n+1}}\) 
is homotopic to the $p^{th}$-power map. Thus 
$p\cdot\pi_{\ast}(S^{2n+1})$ factors through $\pi_{\ast}(S^{2n-1})$. 
Inducting on $n$ therefore shows that $p^{n}\cdot\pi_{m}(S^{2n+1})=0$ 
for all $m>2n+1$. 

The map $\varphi$ was obtained as a consequence of very thorough 
analysis of the homotopy theory of Moore spaces. The mod-$p^{r}$ 
Moore space $P^{m}(p^{r})$ is the homotopy cofiber of the degree $p^{r}$ 
map on~$S^{m-1}$. Let $F$ be the homotopy fiber of the pinch map   
\(q\colon\namedright{P^{m}(p^{r})}{}{S^{m}}\) 
onto the top cell. When $m=2n+1$, Cohen, Moore, and Neisendorfer 
studied the homotopy fibration sequence 
\[\namedddright{\Omega^{2} S^{2n+1}}{}{\Omega F}{} 
    {\Omega P^{2n+1}(p^{r})}{\Omega q}{\Omega S^{2n+1}}\] 
and painstakingly determined a decomposition of $\Omega F$. 
In particular, one of the factors of $\Omega F$ is~$S^{2n-1}$, giving 
a composite  
\(\varphi_{r}\colon\nameddright{\Omega^{2} S^{2n+1}}{} 
    {\Omega F^{2n+1}}{}{S^{2n-1}}\).  
The $r=1$ case is the map~$\varphi$ above. 

It was later shown that the map $\varphi_{r}$ fits in a homotopy 
fibration sequence  
\begin{equation} 
  \label{anickfib} 
  \namedddright{\Omega^{2} S^{2n+1}}{\varphi_{r}}{S^{2n-1}}{} 
    {\anick}{}{\Omega S^{2n+1}}. 
\end{equation} 
This was proved by Anick~\cite{A} for $p\geq 5$ by using Cohen, 
Moore, and Neisendorfer's work as the base-case in a complex 
induction. A much simpler proof was given in~\cite{GT} for $p\geq 3$ 
by using new methods. 

These methods raised the possibility of producing a different  
proof of the exponent result for spheres which depends more 
on properties of spheres than Moore spaces. The purpose of this 
paper is to give such a proof. The idea is to first construct a 
homotopy fibration as in~(\ref{anickfib}) and then obtain the 
exponent result as a consequence. Explicitly, from first principles, 
we prove the following.  
 
\begin{theorem} 
   \label{main} 
   Let $p\geq 3$ and $r\geq 1$. There is a homotopy fibration sequence 
   \[\namedddright{\Omega^{2} S^{2n+1}}{\varphi_{r}}{S^{2n-1}}{} 
     {\anick}{}{\Omega S^{2n+1}}\] 
   with the property that the composition 
   \(\nameddright{\Omega^{2} S^{2n+1}}{\varphi_{r}}{S^{2n-1}}{E^{2}} 
        {\Omega^{2} S^{2n+1}}\) 
   is homotopic to the $p^{r}$-power map. 
\end{theorem} 

Consequently, inducting on the $r=1$ case gives the exponent result. 

\begin{corollary} 
   \label{exp} 
   The $p$-torsion in $\pi_{\ast}(S^{2n+1})$ is annihilated 
   by $p^{n}$.~$\qqed$ 
\end{corollary} 

One advantage of our work compared to Cohen, Moore, and 
Neisendorfer's is that it bypasses all the hard work in~\cite{CMN1} 
to find and then keep track of higher order torsion in the homotopy 
groups of Moore spaces. A second advantage is that it works 
simultaneously for all odd primes, whereas the $p=3$ case 
in~\cite{N1} required technical modifications to the $p\geq 5$ cases 
in~\cite{CMN1,CMN2}. A third advantage is that it constructs Anick's 
fibration at the same time as proving the exponent result. On the 
other hand, by intent, our work says little about the homotopy 
theory of Moore spaces whereas~\cite{CMN1} gives a very thorough  
account of it. 

This work is partly motivated for two other reasons as well. The  
first concerns the homotopy fiber $W_{n}$ of the double suspension. 
A long-standing conjecture is that $W_{n}$ is a double loop space. 
In particular, potentially $W_{n}\simeq\Omega^{2} T^{2np+1}(p)$. 
Such a homotopy equivalence would have deep implications, one of 
which being a determination of many differentials in the $E_{2}$-term  
of the $EHP$-spectral sequence calculating the homotopy groups of 
spheres. The sticking point is that $W_{n}$ naturally arises from an 
$EHP$-perspective while $T^{2np+1}(p)$ naturally arises from a Cohen, 
Moore, Neisendorfer (CMN) perspective, and it is not known how 
to reconcile the two. Our starting point in constructing the 
space \anick\ is to use Gray's~\cite{G} construction of a classifying 
space $BW_{n}$ of $W_{n}$. So this paper can be thought of as an 
attempt to fuse the CMN and $EHP$-perspectives, in the hope that it 
will tell us more about whether $W_{n}\simeq\Omega^{2} T^{2np+1}(p)$. 
The second motivation concerns whether there are analogues of 
Anick's fibration at the prime $2$. For Hopf invariant one reasons 
the space $T^{2n+1}(2)$ cannot exist if $n\notin\{2,4,8\}$, but 
Cohen has conjectured that $T^{2n+1}(2^{r})$ does exist if $r>1$. 
The approach in~\cite{GT} would come close to constructing such 
spaces if it were known that $BW_{n}$ is a homotopy associative 
$H$-space when localized at $2$. However, this is not known and 
is probably not the case. The approach to constructing the space 
\anick\ here when $p$ is odd is a variant of that in~\cite{GT}, 
which bypasses the issue of whether $BW_{n}$ is an $H$-space, and 
so raises the possibility of positively answering Cohen's 
conjecture. This will be pursued in later work.

\section{An outline of the proof of Theorem~\ref{main}} 
\label{sec:proof} 

In this section we prove Theorem~\ref{main}, deferring 
details to later sections in order to make the thrust of 
the proof clear. We begin by defining some spaces and maps, and 
stating their relevant properties. 

Let $p$ be an odd prime and $r\in\mathbb{N}$. Let 
\(p^{r}\colon\namedright{S^{2n+1}}{}{S^{2n+1}}\) 
be the map of degree~$p^{r}$. Let $\curly$ be its homotopy fiber. 
As $p$ is an odd prime, $S^{2n+1}$ is an $H$-space, its $p^{r}$-power 
map is homotopic to the degree~$p^{r}$ map, and the $p^{r}$-power 
map is an $H$-map. In particular, \curly\ is an $H$-space. 
In~\cite{N3} it was shown that the $p^{r}$-power map on \curly\ 
is null homotopic.  
 
In~\cite{G} it was shown that the fiber $W_{n}$ of the double suspension 
has a classifying space~$BW_{n}$ and there are homotopy fibrations 
\[\nameddright{S^{2n-1}}{E^{2}}{\Omega^{2} S^{2n+1}}{\nu}{BW_{n}}\] 
\[\nameddright{BW_{n}}{j}{\Omega^{2} S^{2np+1}}{\phi}{S^{2np-1}}\] 
with the property that the composite $j\circ\nu$ is homotopic 
to $\Omega H$, where $H$ is the $p^{th}$-James-Hopf invariant. Further, 
$BW_{n}$ is an $H$-space and the maps $\nu$ and $j$ are $H$-maps  
(the $p=3$ case of this being proved in~\cite{T2}). 
Based on the fact that $\Omega H$ has order $p$, it was shown 
in~\cite{G} that the composite 
\(\nameddright{BW_{n}}{j}{\Omega^{2} S^{2np+1}}{p} 
    {\Omega^{2} S^{2np+1}}\) 
is null homotopic, implying that $j$ lifts to a map 
\(\overline{S}\colon\namedright{BW_{n}}{}{\Omega^{2} S^{2np+1}\{p\}}\). 
In~\cite{T1}, based on a result in~\cite{S}, it was shown 
that $\overline{S}$ can be chosen to be an $H$-map. In 
Lemma~\ref{BWnfib} we will show that the homotopy fiber of $\overline{S}$ 
is $\Omega W_{np}$. Let $S$ be the composite 
\[S\colon\nameddright{\Omega^{2} S^{2n+1}}{\nu}{BW_{n}}{\overline{S}} 
      {\Omega^{2} S^{2np+1}\{p\}}.\] 
Note that $S$ is an $H$-map as it the composite of $H$-maps. 

Now we turn to setting up the proof of Theorem~\ref{main}. Define 
the space $Y$ and the maps $f$, $g$, and~$h$ by the homotopy 
pullback diagram 
\begin{equation} 
  \label{Ydgrm} 
  \diagram 
    S^{2n-1}\rto^-{f}\ddouble & Y\rto^-{g}\dto^{h} 
       & \Omega W_{np}\dto \\ 
    S^{2n-1}\rto^-{E^{2}} & \Omega^{2} S^{2n+1}\rto^-{\nu}\dto^{S}  
       & BW_{n}\dto^{\overline{S}} \\ 
    & \Omega^{2} S^{2np+1}\{p\}\rdouble & \Omega^{2} S^{2np+1}\{p\}. 
  \enddiagram 
\end{equation}  
As it stands, Cohen, Moore, and Neisendorfer's work tells us that 
there is a lift 
\[\diagram 
      & \Omega^{2} S^{2n+1}\dto^{p^{r}}\dlto \\ 
      S^{2n-1}\rto^-{E^{2}} & \Omega^{2} S^{2n+1} 
  \enddiagram\] 
for each $r\geq 1$. To reproduce this in our case, we will find 
successive lifts 
\begin{equation} 
  \label{twolifts} 
  \diagram 
      & & \Omega^{2} S^{2n+1}\ddto^{p^{r}} 
            \ddlto^{\ell_{1}}\ddllto_{\ell_{2}} \\ 
      & & \\ 
      S^{2n-1}\rto^-{f} & Y\rto^-{h} & \Omega^{2} S^{2n+1}. 
  \enddiagram 
\end{equation} 
Both $\ell_{1}$ and $\ell_{2}$ are constructed as consequences of certain 
extensions. To describe how $\ell_{1}$ comes about, consider the pinch map 
\(\namedright{P^{2n+1}(p^{r})}{q}{S^{2n+1}}\) 
onto the top cell. This factors as the composite 
\(\nameddright{P^{2n+1}(p^{r})}{i}{\curly}{}{S^{2n+1}}\) 
where $i$ is the inclusion of the bottom Moore space. The 
factorization determines an extended homotopy pullback diagram 
\begin{equation} 
  \label{EFdgrm} 
  \diagram 
      \Omega^{2} S^{2n+1}\rto^-{\partial_{E}}\dto 
         & E\rto\ddouble & F\rto\dto & \Omega S^{2n+1}\dto \\ 
      \Omega\curly\rto & E\rto & P^{2n+1}(p^{r})\rto^-{i}\dto^{q} 
         & \curly\dto \\ 
      & & S^{2n+1}\rdouble & S^{2n+1} 
  \enddiagram 
\end{equation} 
which defines the spaces $E$ and $F$ and the map $\partial_{E}$. 

\begin{proposition} 
   \label{Eextension} 
   There is an extension 
   \[\diagram 
         \Omega^{2} S^{2n+1}\rto^-{\partial_{E}}\dto^{S} 
             & E\dldashed|>\tip^-{e_{1}} \\ 
         \Omega^{2} S^{2np+1}\{p\}  & 
     \enddiagram\] 
   for some map $e_{1}$. 
\end{proposition} 

By~(\ref{EFdgrm}) the map $\partial_{E}$ factors through 
$\Omega\curly$. Using this, let $e^{\prime}_{1}$ be the composite 
\[e^{\prime}_{1}\colon\nameddright{\Omega\curly}{}{E}{e_{1}} 
     {\Omega^{2} S^{2np+1}\{p\}}.\]  
Then Proposition~\ref{Eextension} implies that there is a homotopy 
commutative square 
\[\diagram 
     \Omega^{2} S^{2n+1}\rto\dto^{S} 
        & \Omega\curly\dto^-{e^{\prime}_{1}} \\ 
     \Omega^{2} S^{2np+1}\{p\}\rdouble & \Omega^{2} S^{2np+1}\{p\}. 
  \enddiagram\]  
This square determines an extended homotopy pullback diagram 
\begin{equation} 
  \label{l1dgrm} 
  \diagram 
     \Omega^{2} S^{2n+1}\rto^-{\ell_{1}}\ddouble & Y\rto\dto^{h}  
        & X\rto^-{q_{X}}\dto & \Omega S^{2n+1}\ddouble \\ 
     \Omega^{2} S^{2n+1}\rto^-{p^{r}} & \Omega^{2} S^{2n+1}\rto\dto^{S} 
        & \Omega\curly\rto\dto^{e^{\prime}_{1}} 
        & \Omega S^{2n+1} \\ 
     & \Omega^{2} S^{2np+1}\{p\}\rdouble & \Omega^{2} S^{2np+1}\{p\} &   
  \enddiagram 
\end{equation} 
which defines the space $X$ and the maps $q_{X}$ and $\ell_{1}$. 
In particular, $\ell_{1}$ is a choice of a lift of the $p^{r}$-power 
map with the 
additional property that it is a connecting map of a homotopy fibration. 

Further, the factorization of $e^{\prime}_{1}$ through $E$ determines 
a homotopy pullback diagram 
\begin{equation} 
  \label{OmegaPXdgrm} 
  \diagram 
      \Omega\moore\rto^-{i_{X}}\ddouble & X\rto\dto & R\dto \\ 
      \Omega\moore\rto^-{\Omega i}  
         & \Omega\curly\rto\dto^-{e^{\prime}_{1}} & E\dto^{e_{1}} \\ 
      & \Omega^{2} S^{2np+1}\{p\}\rdouble & \Omega^{2} S^{2np+1}\{p\}   
  \enddiagram 
\end{equation}  
which defines the space $R$ and the map $i_{X}$. Note 
that~(\ref{l1dgrm}) and~(\ref{OmegaPXdgrm}) imply that the composite 
\(\nameddright{\Omega\moore}{i_{X}}{X}{q_{X}}{\Omega S^{2n+1}}\) 
is homotopic to~$\Omega q$. The map $i_{X}$ and the homotopy 
$\Omega q\simeq q_{X}\circ i_{X}$ will play an important 
role in proving the existence of the second extension. 

Consider the homotopy fibration 
\(\nameddright{Y}{}{X}{q_{X}}{\Omega S^{2n+1}}\) 
and the map 
\(\namedright{Y}{g}{\Omega W_{np}}\) 
from~(\ref{Ydgrm}). 

\begin{proposition} 
   \label{Xextension}  
   There is an extension 
   \[\diagram  
         Y\rto\dto^{g} & X\dldashed|>\tip^-{e_{2}} \\
         \Omega W_{np} &
     \enddiagram\]  
   for some map $e_{2}$. 
\end{proposition} 

Proposition~\ref{Xextension} implies that there is an extended 
homotopy pullback diagram    
\begin{equation} 
  \label{Tpb} 
  \diagram 
     \Omega^{2} S^{2n+1}\rto^-{\ell_{2}}\ddouble 
        & S^{2n-1}\rto\dto^{f} & \anick\rto\dto 
        & \Omega S^{2n+1}\ddouble \\ 
     \Omega^{2} S^{2n+1}\rto^-{\ell_{1}} & Y\rto\dto^{g} 
        & X\rto^-{q_{X}}\dto & \Omega S^{2n+1} \\ 
     & \Omega W_{np}\rdouble & \Omega W_{np} & 
  \enddiagram 
\end{equation}  
which defines the space \anick\ and the map $\ell_{2}$. 

\medskip\noindent 
\textit{Proof of Theorem~\ref{main}}: 
The top row in~(\ref{Tpb}) is the asserted homotopy fibration. 
Relabel $\ell_{2}$ as $\varphi_{r}$. All that is left to check is 
that $\varphi_{r}$ has the property that $E^{2}\circ\varphi_{r}$ 
is homotopic to the $p^{r}$-power map on $\Omega^{2} S^{2n+1}$. 
This follows from juxtaposing the squares of connecting maps 
in~(\ref{l1dgrm}) and~(\ref{Tpb}), giving a homotopy commutative 
diagram 
\[\diagram 
     \Omega^{2} S^{2n+1}\rto^-{\varphi_{r}}\ddouble & S^{2n-1}\dto^{f} \\ 
     \Omega^{2} S^{2n+1}\rto^-{\ell_{1}}\ddouble & Y\dto^{h} \\ 
     \Omega^{2} S^{2n+1}\rto^-{p^{r}} & \Omega^{2} S^{2n+1},   
  \enddiagram\] 
and noting that by~(\ref{Ydgrm}) the composite $h\circ f$ is homotopic 
to $E^{2}$. 
$\qqed$ 
\medskip 

We close this section with some remarks regarding the proof  
of Theorem~\ref{main}. The key is establishing the existence of the 
extensions in Propositions~\ref{Eextension} and~\ref{Xextension}. 
These are obtained by filtering $E$ and $X$ respectively by certain 
homotopy pushouts. Inductively, extensions are produced  
as pushout maps. To determine that a pushout map exists, we use 
a slight generalization of a theorem in~\cite{GT}. One of the 
hypotheses of this theorem involves a certain map in the pushout 
being divisible by $p^{r}$, or having its suspension divisible 
by $p^{r}$. This is played off against the fact that the $p^{th}$-power 
maps on $\Omega^{2} S^{2np+1}\{p\}$ and $\Omega W_{np}$ are 
null homotopic in order to annihilate the obstructions to the 
existence of a pushout map.  

It has already been pointed out that $S^{2np+1}\{p\}$ has a 
null homotopic $p^{th}$-power map, implying the same for 
$\Omega^{2} S^{2np+1}\{p\}$. The case of $\Omega W_{np}$ needs 
to be stated clearly. In~\cite{CMN1} it was shown that 
the $p^{th}$-power map on $W_{n}$ is null homotopic, but this 
used the fact that the $p^{th}$-power map on $\Omega^{2} S^{2n+1}$ 
factors as the composite 
\(\nameddright{\Omega^{2} S^{2n+1}}{\varphi}{S^{2n-1}}{E^{2}} 
    {\Omega^{2} S^{2n+1}}\). 
In~\cite{T2} the same factorization was used to show the stronger result 
that the $p^{th}$-power map on $BW_{n}$ is null homotopic. As we are 
trying to give a new proof of this factorization, we have to suppress 
knowledge of these exponent results. However, as will be seen in 
Section~\ref{sec:background}, the method of proof in~\cite{T2} 
can be applied to an alternate factorization of the $p^{th}$-power 
map on $\Omega^{2} S^{2np+1}$ as the composite 
\(\nameddright{\Omega^{2} S^{2np+1}}{\phi}{S^{2np-1}}{E^{2}} 
    {\Omega^{2} S^{2np+1}}\), 
which is independent of~\cite{CMN1}, in order to show that the 
$p^{th}$-power map on $BW_{np}$ is null homotopic. 

The remainder of the paper is organized as follows.  
Section~\ref{sec:background} gives some preliminary background results 
related to $BW_{n}$. In Section~\ref{sec:extensionthm} the generalized 
version of the extension theorem in~\cite{GT} is given. The 
$p^{r}$-divisibility condition is established in two contexts in 
Sections~\ref{sec:pdiv1} and~\ref{sec:pdiv2}, corresponding to 
the contexts of Propositions~\ref{Eextension} and~\ref{Xextension}. 
The two propositions are then proved in Section~\ref{sec:extensions}. 
Finally, in Section~\ref{sec:consequences} we give some concluding 
remarks and consequences.

\section{Background results} 
\label{sec:background} 

In this section we establish some background material related 
to $BW_{n}$. First, consider the map  
\(\namedright{\Omega^{2} S^{2np+1}}{\phi}{S^{2np-1}}\) 
whose homotopy fiber is $BW_{n}$. In~\cite{T1} it was shown that 
the composite
\(\nameddright{\Omega^{2} S^{2np+1}}{\phi}{S^{2np-1}}{E^{2}}
     {\Omega^{2} S^{2np+1}}\)
is homotopic to the $p^{th}$-power map. The following lemma 
is an adaptation of~\cite[1.2]{T2}. 

\begin{lemma}
   \label{BWnpexp}
   The $p^{th}$-power map on $BW_{np}$ is null homotopic. 
\end{lemma}

\begin{proof} 
Since $BW_{np}$ is an $H$-space, it suffices to show that 
\(\namedright{\Sigma BW_{np}}{\Sigma p}{\Sigma BW_{np}}\) 
is null homotopic. As the $p^{th}$-power map on $\Omega^{2} S^{2np+1}$ 
factors as    
\(\nameddright{\Omega^{2} S^{2np+1}}{\phi}{S^{2np-1}}{E^{2}} 
   {\Omega^{2} S^{2np+1}}\), 
the composite $\nu\circ p$ factors as $\nu\circ E^{2}\circ\phi$. 
Since $\nu$ and $E^{2}$ are consecutive maps in a homotopy fibration, 
their composite is null homotopic and so $\nu\circ p$ is null homotopic. 
Since $\nu$ is an $H$-map, $p\circ\nu\simeq\nu\circ p$, and so 
$p\circ\nu$ is null homotopic. Let $C$ be the homotopy cofiber 
of $\nu$. Then the null homotopy for $p\circ\nu$ implies that there 
is an extension 
\[\diagram 
      \Omega^{2} S^{2n+1}\rto^-{\nu} & BW_{n}\rto^-{t}\dto^{p} 
         & C\dldashed|>\tip^-{\lambda} \\ 
      & BW_{n} & 
  \enddiagram\] 
for some map $\lambda$. By~\cite{G}, the map $\Sigma^{2}\nu$ has 
a right homotopy inverse. This implies that there is a homotopy 
decomposition  
$\Sigma^{2}\Omega^{2} S^{2n+1}\simeq\Sigma C\vee\Sigma^{2} BW_{n}$.   
Thus $\Sigma t$ is null homotopic and so $\Sigma p$ is null homotopic.
\end{proof} 

Next, consider the $H$-map 
\(\namedright{BW_{n}}{\overline{S}}{\Omega^{2} S^{2np+1}\{p\}}\) 
which lifts 
\(\namedright{BW_{n}}{j}{\Omega^{2} S^{2np+1}}\). 
In Section~\ref{sec:proof} it was stated that the homotopy 
fiber of $\overline{S}$ is $\Omega W_{np}$. We prove this now.     

\begin{lemma} 
   \label{BWnfib} 
   There is a homotopy fibration 
   \(\nameddright{\Omega W_{np}}{}{BW_{n}}{\overline{S}} 
       {\Omega^{2} S^{2np+1}\{p\}}\). 
\end{lemma} 

\begin{proof} 
Consider the homotopy pullback diagram 
\[\diagram 
      M\rto\ddouble & \Omega S^{2np-1}\rto^-{s}\dto 
         & \Omega^{3} S^{2np+1}\dto \\ 
      M\rto & BW_{n}\rto^-{\overline{S}}\dto^{j} 
         & \Omega^{2} S^{2np+1}\{p\}\dto \\ 
      & \Omega^{2} S^{2np+1}\rdouble & \Omega^{2} S^{2np+1} 
  \enddiagram\] 
which defines the space $M$ and the map $s$. Since this 
is a pullback of $H$-spaces and $H$-maps, the induced map $s$ 
is an $H$-map. The James construction~\cite{J} implies that any $H$-map 
\(t\colon\namedright{\Omega\Sigma X}{}{Z}\), 
where $Z$ is homotopy associative, is determined uniquely, 
up to homotopy, by $t\circ E$. In our case, $s\circ E$ is 
of degree~$1$ in homology and so $s\circ E\simeq E^{3}$. 
On the other hand, 
\(\namedright{\Omega S^{2np-1}}{\Omega E}{\Omega^{3} S^{2np+1}}\) 
is another $H$-map with the property that $\Omega E\circ E\simeq E^{3}$. 
Thus $s\simeq\Omega E$. Hence $M\simeq\Omega W_{np}$. 
\end{proof} 

Finally, we need to show the existence of a certain homotopy 
commutative diagram describing a property of the $H$-map 
\(\namedright{\Omega^{2} S^{2n+1}}{S}{\Omega^{2} S^{2np+1}\{p\}}\).  
In general, if $Z$ is an $H$-space and there are maps 
\(f\colon\namedright{X}{}{Z}\) 
and 
\(g\colon\namedright{Y}{}{Z}\), 
let $f\cdot g$ be the product of $f$ and $g$, obtained as the 
composite 
\[f\cdot g\colon\nameddright{X\times Y}{f\times g}{Z\times Z}{m}{Z}\] 
where $m$ is the multiplication on $Z$.  

\begin{lemma} 
   \label{Ezerocase} 
   There is a homotopy commutative diagram 
   \[\diagram 
        S^{2n-1}\times\Omega^{2} S^{2n+1} 
              \rto^-{E^{2}\cdot 1}\dto^{\pi_{2}} 
          & \Omega^{2} S^{2n+1}\dto^{S} \\ 
        \Omega^{2} S^{2n+1}\rto^-{S} & \Omega^{2} S^{2n+1}\{p\} 
     \enddiagram\] 
   where $\pi_{2}$ is the projection. 
\end{lemma} 

\begin{proof} 
Observe that the composite 
\(\nameddright{S^{2n-1}}{E^{2}}{\Omega^{2} S^{2n+1}}{S} 
     {\Omega^{2} S^{2np+1}\{p\}}\) 
is null homotopic by connectivity. The lemma now follows 
since $S$ is an $H$-map. 
\end{proof}

\section{An extension theorem} 
\label{sec:extensionthm} 

Let 
\(\nameddright{A}{}{E^{\prime}}{}{E}\) 
be a homotopy cofibration and suppose there is a map 
\(\namedright{F}{}{E}\). 
Define spaces $F^{\prime}$, $Q$, and $X$ by the diagram of 
iterated homotopy pullbacks 
\begin{equation} 
  \label{DLdgrm} 
  \diagram 
     X\rdouble\dto & X\rdouble\dto & X\dto \\ 
     Q\rto\dto & F^{\prime}\rto\dto & F\dto \\ 
     A\rto & E^{\prime}\rto & E. 
  \enddiagram 
\end{equation}  
The following Lemma can be found in~\cite{G}.   

\begin{lemma} 
   \label{DLpo} 
   There is a homotopy equivalence  
   \(\namedright{Q}{\simeq}{A\times X}\) 
   such that the map 
   \(\namedright{Q}{}{A}\) 
   becomes the projection, and there is a homotopy pushout 
   \[\diagram 
        Q\simeq A\times X\rto\dto^{\pi_{2}} & F^{\prime}\dto \\ 
        X\rto & F. 
     \enddiagram\]  
\end{lemma} 

\begin{remark} 
\label{DLpoaction} 
Note that there may be many choices of the homotopy equivalence 
for $Q$, and therefore many possible homotopy classes for the map 
\(\namedright{A\times X}{}{F^{\prime}}\). 
One situation where additional control can be imposed is in the 
case of a principal fibration. Suppose $X=\Omega X^{\prime}$ and 
there is a homotopy fibration sequence 
\(\namedddright{\Omega X^{\prime}}{}{F^{\prime}}{}{E^{\prime}} 
     {}{X^{\prime}}\). 
Then there is a homotopy action 
\(\theta\colon\namedright{F^{\prime}\times\Omega X^{\prime}}{} 
   {F^{\prime}}\). 
As observed in~\cite{GT}, the decomposition $Q\simeq A\times X$ can 
then be chosen so that the map 
\(\namedright{A\times X}{}{F^{\prime}}\) 
is homotopic to the composite 
\(\lnameddright{A\times X}{a\times 1}{F^{\prime}\times X}{\theta}
   {F^{\prime}}\), 
where $a$ is any choice of a lift of 
\(\namedright{A}{}{E^{\prime}}\) 
to $F^{\prime}$. Note that such a lift exists by the decomposition 
of $Q$, and there may be many choices of $a$.  
\end{remark} 

Theorem~\ref{extension}, a slight generalization 
of~\cite[2.3]{GT}, gives conditions for when the homotopy pushout 
in Lemma~\ref{DLpo} can be used to extend a map 
\(\namedright{F^{\prime}}{}{Z}\) 
to a map 
\(\namedright{F}{}{Z}\).  
As notation, if $A$ is a co-$H$ space, let 
\(\underline{p}^{r}\colon\namedright{A}{}{A}\) 
be the map of degree~$p^{r}$. 

\begin{theorem} 
   \label{extension} 
   Given a diagram as in~(\ref{DLdgrm}) in which $A$ is a suspension. 
   Suppose the homotopy equivalence $Q\simeq A\times X$ in 
   Lemma~\ref{DLpo} can be chosen so that the map 
   \(\namedright{Q\simeq A\times X}{}{F^{\prime}}\) 
   is homotopic to a composite 
   \[\nameddright{A\times X}{a\times 1}{M\times X}{}{F^{\prime}}\] 
   for some space $M$, where $a$ has the property that 
   $\Sigma a\simeq t\circ\underline{p}^{r}$ for some map $t$. 
   Let $Z$ be a homotopy associative $H$-space whose $p^{r}$-power map 
   is null homotopic and suppose there is a map  
   \(\namedright{F^{\prime}}{}{Z}\). 
   Then there is an extension 
   \[\diagram 
        F^{\prime}\rto\dto & Z\ddouble \\  
        F\rto & Z. 
     \enddiagram\]  
\end{theorem} 

\begin{proof} 
The theorem was proved in~\cite[2.3]{GT} in the special case where 
$X=\Omega X^{\prime}$ and there is a homotopy fibration sequence 
\(\namedddright{\Omega X^{\prime}}{}{F^{\prime}}{f^{\prime}} 
   {E^{\prime}}{}{X^{\prime}}\). 
By Remark~\ref{DLpoaction}, the decomposition 
$Q\simeq A\times X$ may be chosen so that the map 
\(\namedright{Q\simeq A\times X}{}{F^{\prime}}\) 
is homotopic to the composite 
\[\overline{\theta}\colon\nameddright{A\times X}{a\times 1} 
   {F^{\prime}\times X}{}{F^{\prime}},\] 
where $a$ is some lift of 
\(\namedright{A}{}{E}\).  
So in this case, $M=F^{\prime}$. The existence of the extension 
is proved in~\cite[2.6]{GT} given the composite $\overline{\theta}$, 
the fact that $\Sigma a$ is divisible by $p^{r}$, and the $H$-space 
properties of $Z$. 

In the more general case, the proof of~\cite[2.6]{GT} holds without 
change once $\overline{\theta}$ is replaced by the given composite 
\(\nameddright{A\times X}{a\times 1}{M\times X}{}{F^{\prime}}\).   
\end{proof}

\section{$p^{r}$-divisibility I} 
\label{sec:pdiv1} 

In the next two sections we prove $p^{r}$-divisibility results 
which will eventually be used as input into  
Theorem~\ref{extension}. Assume homology is taken with mod-$p$ coefficients. 

The first $p^{r}$-divisibility result concerns the spaces $E$ and $F$    
in the homotopy pullback diagram   
\begin{equation} 
  \label{EFdgrm2} 
  \diagram 
     E\rto\ddouble & F\rto\dto & \Omega S^{2n+1}\dto \\ 
     E\rto & P^{2n+1}(p^{r})\rto^-{i}\dto^{q} & \curly\dto \\ 
     & S^{2n+1}\rdouble & S^{2n+1}  
  \enddiagram 
\end{equation} 
introduced in Section~\ref{sec:proof}.  
A straightforward calculation of the mod-$p$ homology 
Serre spectral sequence for the homotopy fibration 
\(\nameddright{\Omega S^{2n+1}}{}{F}{}{P^{2n+1}(p^{r})}\) 
shows that $\rhlgy{F}$ is generated as a vector space by 
elements $x_{2nk}$ in degrees $2nk$ for $k\geq 1$. Thus $F$ has 
a natural filtration by skeleta. Let $F_{k}$ be the $2nk$-skeleton 
of $F$. Then for each $k\geq 1$ there is a homotopy cofibration 
\[\lnameddright{S^{2nk-1}}{g_{k-1}}{F_{k-1}}{}{F_{k}}.\]  
The filtration on $F$ lets us put a filtration on $E$. Define 
the spaces $E_{k-1}$ and $D_{k-1}$ by the homotopy pullback diagram 
\begin{equation} 
  \label{EkFk} 
  \diagram 
      D_{k-1}\rto\ddouble & E_{k-1}\rto\dto & E\dto \\ 
      D_{k-1}\rto & F_{k-1}\rto\dto & F\dto \\ 
      & \Omega S^{2n+1}\rdouble & \Omega S^{2n+1}. 
  \enddiagram 
\end{equation} 
Since 
\(\namedright{F_{k-1}}{}{F}\) 
is a skeletal inclusion, a Serre spectral sequence calculation 
immediately implies that $D_{k-1}$ is $(2nk-2)$-connected and has 
a single cell in dimension $2nk-1$. Further, the composite 
\(S^{2nk-1}\hookrightarrow\namedright{D_{k-1}}{}{F_{k-1}}\) 
is homotopic to~$g_{k-1}$. Let $\overline{g}_{k-1}$ be the composite 
\[\overline{g}_{k-1}\colon S^{2nk-1}\hookrightarrow\namedright{D_{k-1}} 
     {}{E_{k-1}}.\] 
In Proposition~\ref{Edivbyp} we will show that $\overline{g}_{k-1}$ 
is divisible by $p^{r}$, given mild conditions.   

\begin{proposition} 
   \label{Edivbyp} 
   If $k>1$ then there is a homotopy commutative diagram 
   \[\diagram 
         S^{2nk-1}\rto^-{\overline{g}_{k-1}}\dto^{\underline{p}^{r}} 
             & E_{k-1}\ddouble \\ 
         S^{2nk-1}\rto^-{\overline{h}_{k-1}} & E_{k-1} 
     \enddiagram\] 
   for some map $\overline{h}_{k-1}$.  
\end{proposition}  

Proposition~\ref{Edivbyp} appeared in~\cite{GT} and relied heavily on 
the work of Cohen, Moore, and Neisendorfer~\cite{CMN1}. However, we will 
give a different proof which is self-contained. While some of the basic 
tools will be the same, we will require much less homological algebra 
and virtually none of the differential graded Lie algebra 
techniques from~\cite{CMN1} which underpinned the proof in~\cite{GT}. 

We begin with a small amount of mod-$p$ homotopy theory; the primary  
reference is~\cite{N1}. In general, suppose we are given maps 
\(f\colon\namedright{P^{s}(p^{r})}{}{\Omega X}\) 
and 
\(g\colon\namedright{P^{t}(p^{r})}{}{\Omega X}\). 
As the smash $P^{s}(p^{r})\wedge P^{t}(p^{r})$ 
is homotopy equivalent to the wedge $P^{s+t}(p^{r})\vee P^{s+t-1}(p^{r})$,  
we can define the mod-$p^{r}$ Samelson product of $f$ and $g$ 
as the composite 
\[\langle f,g\rangle\colon\nameddright{P^{s+t}(p^{r})}{} 
   {P^{s}(p^{r})\wedge P^{t}(p^{r})}{[f,g]}{\Omega X}\] 
where $[f,g]$ is the ordinary Samelson product of $f$ and $g$. 
In our case, let 
\(\nu\colon\namedright{P^{2n}(p^{r})}{}{\Omega P^{2n+1}(p^{r})}\) 
be the adjoint of the identity map and let 
\(\mu\colon\namedright{S^{2n-1}}{}{\Omega P^{2n+1}(p^{r})}\) 
be the inclusion of the bottom cell. Let $\mathfrak{ad}^{0}=\mu$, 
$\mathfrak{ad}^{1}$ be the ordinary Samelson product $[\nu,\mu]$, 
and for $j>1$ let  
\(\mathfrak{ad}^{j-1}:\namedright{P^{2nj-1}(p^{r})}{}{\Omega P^{2n+1}}\) 
be the mod-$p^{r}$ Samelson product defined inductively by 
$\mathfrak{ad}^{j-1}=\langle\nu,\mathfrak{ad}^{j-2}\rangle$. 

Homologically, since $P^{2n+1}(p^{r})$ is the suspension 
of $P^{2n}(p^{r})$, the Bott-Samelson theorem implies that 
$\hlgy{\Omega P^{2n+1}(p^{r})}\cong T(\rhlgy{P^{2n}(p^{r})})$, 
where $T(\ )$ denotes the free tensor algebra. A basis for 
$\rhlgy{P^{2n}(p^{r})}$ is given by $\{u,v\}$ where $u$ and 
$v$ are in degrees $2n-1$ and $2n$ respectively. As a matter of 
book-keeping, recall the well known fact that if $V$ is a graded 
module over a field then $T(V)\cong UL\langle V\rangle$, where 
the right side is the universal enveloping algebra of the free 
Lie algebra generated by $V$. In our case, we have 
$\hlgy{\Omega P^{2n+1}(p^{r})}\cong UL\langle u,v\rangle$. 
Picking out particular elements in $L\langle u,v\rangle$, let 
$ad^{0}=u$ and for $j\geq 1$ define $ad^{j}$ inductively by 
$ad^{j}=[v,ad^{j-1}]$.   

In general, given a map 
\(f:\namedright{P^{s}(p^{r})}{}{\Omega X}\), 
the mod-$p^{r}$ Hurewicz homomorphism is defined by sending $f$ 
to $f_{\ast}(\iota)$, where $\iota$ is the generator in 
$H_{s}(P^{s}(p^{r}))$. In our case, the ordinary Hurewicz image of $\mu$ 
is $u$ and the mod-$p^{r}$ Hurewicz image of $\nu$ is $v$. The 
mod-$p^{r}$ Hurewicz image of 
\(\llnamedright{P^{2nj-1}(p^{r})}{\mathfrak{ad}^{j-1}} 
     {\Omega P^{2n+1}(p^{r})}\) 
is therefore $ad^{j-1}$. 

With these ingredients in place, consider the homotopy fibration 
\(\nameddright{\Omega E}{}{\Omega\moore}{\Omega i}{\Omega\curly}\). 
For $j>1$ the map $\mathfrak{ad}^{j-1}$ is defined via a Samelson 
product and so composes trivially with the loop map $\Omega i$ 
since $\Omega\curly$ is homotopy commutative. Thus there is a lift 
\[\diagram 
     & \Omega E\dto \\ 
     P^{2nk-1}(p^{r})\rto^-{\mathfrak{ad}^{j-1}}\urto^-{l_{j-1}} 
        & \Omega\moore 
  \enddiagram\] 
for some map $l_{j-1}$. Let $l^{\prime}_{j-1}$ be the composite 
\(l^{\prime}_{j-1}\colon\nameddright{P^{2nk-1}(p^{r})}{l_{j-1}} 
     {\Omega E}{}{\Omega F}\).  
Then~(\ref{EFdgrm2}) implies that there is a homotopy commutative 
diagram  
\[\diagram 
     & \Omega F\dto \\ 
     P^{2nk-1}(p^{r}) 
           \rto^-{\mathfrak{ad}^{j-1}}\urto^-{l^{\prime}_{j-1}} 
        & \Omega\moore. 
  \enddiagram\] 
Note that for $\Omega F$ there is an additional lift. By 
connectivity, the map 
\(\namedright{S^{2n-1}}{\mathfrak{ad}^{0}}{\Omega\moore}\) 
lifts to a map 
\(l^{\prime}_{0}\colon\namedright{S^{2n-1}}{}{\Omega F}\). 

In particular, the Hurewicz image $ad^{0}$ of $\mathfrak{ad}^{0}$ 
factors through $\hlgy{\Omega F}$, as do the mod-$p^{r}$ Hurewicz 
images $ad^{j-1}$ of $\mathfrak{ad}^{j-1}$ for each $j>1$. 
Since the skeletal inclusion 
\(\namedright{F_{k-1}}{}{F}\) 
is $(2nk-1)$-connected, for dimension and connectivity reasons 
the maps $l^{\prime}_{j-1}$ factor through $\Omega F_{k-1}$ 
for $0\leq j-1<k-1$. Thus the submodule $W_{k-1}=\{ad^{j-1}\}_{j=1}^{k-1}$ 
of Hurewicz images in $\hlgy{\Omega P^{2n+1}(p^{r})}$ 
factors through $\hlgy{\Omega F_{k-1}}$. Extending multiplicatively 
gives a commative diagram 
\[\diagram 
     T(W_{k-1})\rto\dto^{t_{k-1}} & UL\langle u,v\rangle\dto^{\cong} \\ 
     \hlgy{\Omega F_{k-1}}\rto & \hlgy{\Omega P^{2n+1}(p^{r})}  
  \enddiagram\] 
which defines the map $t_{k-1}$. Since $T(W_{k})\cong UL\langle W_{k}\rangle$, 
and $L\langle W_{k}\rangle$ is a sub-Lie algebra of $L\langle u,v\rangle$, 
the Poincar\'{e}-Birkhoff-Witt Theorem implies that 
$UL\langle W_{k}\rangle$ is a subalgebra of $UL\langle u,v\rangle$. 
Thus the upper direction around the diagram is an injection, and 
so $t_{k-1}$ is an injection. In the next Lemma we show that $t_{k-1}$  
is actually an isomorphism. 

\begin{lemma} 
   \label{loopFk-1hlgy} 
   The map 
   \(\namedright{T(W_{k-1})}{t_{k-1}}{\hlgy{\Omega F_{k-1}}}\) 
   is an algebra isomorphism. 
\end{lemma} 

\begin{proof} 
Since $t_{k-1}$ is an algebra injection, to prove the lemma 
it suffices to show that $T(W_{k-1})$ and $\hlgy{\Omega F_{k-1}}$ 
have the same Euler-Poincar\'{e} series. Denote the Euler-Poincar\'{e} 
series of a graded module $A$ by $P(A)$. As $t_{k-1}$ is an injection,  
we have $P(T(W_{k-1}))\leq P(\hlgy{\Omega F_{k-1}})$. 
On the other hand, the cobar construction implies that there is 
an isomorphism 
$\hlgy{\Omega F_{k-1}}\cong\hlgy{T(V_{k-1}),d}$, 
where $V_{k-1}=\Sigma^{-1}\hlgy{F_{k-1}}$ is a desuspension of the module 
$\hlgy{F_{k-1}}$, $d$ is an appropriate differential on the tensor 
algebra, and homology has been taken with respect to that differential. 
In particular, the Euler-Poincar\'{e} series $P(\hlgy{T(V_{k-1}),d})$ 
is bounded above by $P(T(V_{k-1}))$. Next, since $F_{k-1}$ has a single 
cell in dimension $2nj$ for each $1\leq j\leq k-1$, a basis for $V_{k-1}$ 
is $\{x_{2nj-1}\}_{j=1}^{k-1}$, where $x_{2nj-1}$ is in degree $2nj-1$. 
Observe that there is an abstract isomorphism of graded modules 
$W_{k-1}=V_{k-1}$ and so $P(T(V_{k-1}))=P(T(W_{k-1}))$. Putting all 
this together gives inequalities  
\[P(T(W_{k-1}))\leq P(\hlgy{\Omega F_{k-1}}) = P(\hlgy{(T(V_{k-1}),d)} 
    \leq P(T(V_{k-1}))= P(T(W_{k-1}))\]   
and so $P(T(W_{k-1})=P(\hlgy{\Omega F_{k-1}})$, as required.  
\end{proof} 

Inductively, we obtain the following 

\begin{corollary} 
   \label{OmegaFhlgy} 
   There is an isomorphism $\hlgy{\Omega F}\cong T(W_{\infty})$, 
   where $W_{\infty}=\{ad^{j-1}\}_{j=1}^{\infty}$, 
   and the map 
   \(\namedright{\Omega F_{k-1}}{}{\Omega F}\) 
   is modelled homologically by the inclusion of tensor algebras 
   \(\namedright{T(W_{k-1})}{}{T(W_{\infty})}\). 
   $\qqed$ 
\end{corollary} 

Corollary~\ref{OmegaFhlgy} implies that the mod-$p^{r}$ Hurewicz 
image of 
\(\namedright{P^{2nk-1}(p^{r})}{l^{\prime}_{k-1}}{\Omega F}\) 
can be regarded as $ad^{k-1}$. By~(\ref{EkFk}), there is a 
homotopy fibration 
\(\nameddright{\Omega F_{k-1}}{}{\Omega F}{}{D_{k-1}}\). 
Recall that $D_{k-1}$ is $(2nk-2)$-connected and has a single 
cell in dimension $2nk-1$. The following lemma implies that 
in homology $ad^{k-1}$ is mapped to the generator in 
$H_{2nk-1}(D_{k-1})$.   

\begin{lemma} 
   \label{adhurewicz} 
   The composite 
   \(\nameddright{P^{2nk-1}(p^{r})}{l^{\prime}_{k-1}} 
        {\Omega F}{}{D_{k-1}}\) 
   is degree one in $H_{2nk-1}(\ )$. 
\end{lemma} 

\begin{proof} 
By Corollary~\ref{OmegaFhlgy}, the map 
\(\namedright{\Omega F_{k-1}}{}{\Omega F}\) 
is modelled homologically by the inclusion of tensor algebras 
\(\namedright{T(W_{k-1})}{}{T(W_{\infty})}\). 
Since $W_{\infty}=W_{k-1}\oplus V$, where 
$V=\{ad^{j-1}\}_{j=k}^{\infty}$, there is a free product 
decomposition $T(W_{\infty})\cong T(W_{k-1})\coprod T(V)$. 
In particular, this implies that there is a decomposition 
$T(W_{\infty})\cong T(W_{k-1})\otimes N$ 
of left $T(W_{k-1})$-modules for some module $N$.  
Consequently, the Eilenberg-Moore spectral sequence 
for the homotopy fibration 
\(\nameddright{\Omega F_{k-1}}{}{\Omega F_{k}}{}{D_{k-1}}\) 
which converges to $\hlgy{D_{k-1}}$ collapses at $E^{2}$. That is,  
\[E^{2}=\mbox{Tor}^{T(W_{k})}(\zmodp,T(W_{\infty}))= 
    \zmodp\otimes_{T(W_{k})} T(W_{\infty})=N,\]   
and so $\hlgy{D_{k-1}}\cong N$. This implies that the 
element $ad^{k-1}\in V\cap N$ is in $\hlgy{D_{k-1}}$. But 
$ad^{k-1}$ is the mod-$p^{r}$ Hurewicz image of the map  
\(\namedright{P^{2nk-1}(p^{r})}{l^{\prime}_{k-1}}{\Omega F}\), 
and so the composite 
\(\nameddright{P^{2nk-1}(p^{r})}{l^{\prime}_{k-1}} 
     {\Omega F}{}{D_{k-1}}\) 
is degree one in $H_{2nk-1}(\ )$. 
\end{proof}   

We are now ready to prove Proposition~\ref{Edivbyp} 
\medskip 

\noindent 
\textit{Proof of Proposition~\ref{Edivbyp}}: 
We begin with some geometric arguments to reduce the proof to a 
homology calculation. Define the space $M_{k-1}$ by the homotopy 
pullback diagram 
\begin{equation} 
  \label{MNdgrm1} 
  \diagram 
     D_{k-1}\rto\ddouble & M_{k-1}\rto\dto & \Omega\curly\dto \\ 
     D_{k-1}\rto & E_{k-1}\rto\dto & E\dto \\  
     & \moore\rdouble & \moore.  
  \enddiagram 
\end{equation}  
Taking vertical connecting maps and reorienting, we obtain 
a homotopy pullback diagram 
\begin{equation} 
   \label{MNdgrm2} 
   \diagram 
     \Omega E\rto^-{\partial_{1}}\dto & D_{k-1}\rto\dto 
       & E_{k-1}\ddouble \\  
     \Omega\moore\rto\dto^{\Omega i} & M_{k-1}\rto\dto & E_{k-1} \\ 
     \Omega\curly\rdouble & \Omega\curly & 
  \enddiagram 
\end{equation} 
which defines the map $\partial_{1}$. 

By definition, the map 
\(\namedright{P^{2nk-1}(p^{r})}{l_{k-1}}{\Omega E}\) 
is a lift of   
\(\llnamedright{P^{2nk-1}(p^{r})}{\mathfrak{ad}^{j-1}} 
    {\Omega P^{2n+1}(p^{r})}\). 
Let $\lambda^{\prime}$ be the composite 
\[\lambda^{\prime}\colon\nameddright{P^{2nk-1}(p^{r})}{l_{k-1}} 
    {\Omega E}{\partial_{1}}{D_{k-1}}.\] 
Since $D_{k-1}$ is $(2nk-2)$-connected, 
the restriction of $\lambda^{\prime}$ to the $(2nk-2)$-cell of 
$P^{2nk-1}(p^{r})$ is null homotopic, implying that $\lambda^{\prime}$ 
factors as a composite 
\(\nameddright{P^{2nk-1}(p^{r})}{q}{S^{2nk-1}}{\lambda} 
     {D_{k-1}}\), 
for some map~$\lambda$. 

Consider the diagram 
\begin{equation} 
  \label{llambdadgrm} 
  \diagram 
     P^{2nk-1}(p^{r})\rto^-{q}\dto^{l_{k-1}} 
      & S^{2nk-1}\rto^-{\underline{p}^{r}}\dto^{\lambda} 
      & S^{2nk-1}\dto^{h_{k-1}} \\ 
    \Omega E\rto^-{\partial_{1}} & D_{k-1}\rto & E_{k-1}. 
  \enddiagram 
\end{equation}  
The conclusion of the previous paragraph implies that the 
left square homotopy commutes. Since the lower row is a homotopy 
fibration, the homotopy commutativity of the left square implies 
that $\lambda\circ q$ composes trivially into $E_{k-1}$. 
Thus there is an extension along the homotopy cofiber of $q$ 
which makes the right square homotopy commute for some 
map $h_{k-1}$. 

We claim that the composite $\partial_{1}\circ l_{k-1}$ 
is degree one in $H_{2nk-1}(\ )$. Suppose this is the case. Then 
the left square in~(\ref{llambdadgrm}) implies that 
$\lambda$ is degree one in $H_{2nk-1}(\ )$. 
Thus $\lambda$ is homotopic to the inclusion of the bottom 
cell into $D_{k-1}$. So by the definition of $\overline{g}_{k-1}$, 
the composite 
\(\namedddright{S^{2nk-1}}{\lambda}{D_{k-1}}{}{M_{k-1}} 
    {}{E_{k-1}}\) 
is homotopic to $\overline{g}_{k-1}$. Thus the right square 
in~(\ref{llambdadgrm}) implies that 
$\overline{g}_{k-1}\simeq h_{k-1}\circ\underline{p}^{r}$, 
as asserted by the proposition. 

It remains to show that $\partial_{1}\circ l_{k-1}$ is degree one 
in $H_{2nk-1}(\ )$. Refining a bit, by~(\ref{EkFk}) there is an 
extended homotopy pullback diagram 
\[\diagram 
     \Omega E\rto^-{\partial_{1}}\dto & D_{k-1}\rto\ddouble 
        & E_{k-1}\rto\dto & E\dto \\ 
     \Omega F\rto^-{\partial_{2}} & D_{k-1}\rto & F_{k-1}\rto & F 
  \enddiagram\] 
which defines the map $\partial_{2}$. By definition, $l^{\prime}_{k-1}$ 
is the composite 
\(\nameddright{P^{2nk-1}(p^{r})}{l_{j-1}} 
     {\Omega E}{}{\Omega F}\). 
So to show that $\partial_{1}\circ l_{k-1}$ is degree one in 
$H_{2nk-1}(\ )$ it is equivalent to show that 
$\partial_{2}\circ l^{\prime}_{k-1}$ is degree one 
in $H_{2nk-1}(\ )$. But this is the case by Lemma~\ref{adhurewicz}. 
$\qqed$

\section{$p^{r}$-divisibility II} 
\label{sec:pdiv2} 

The second $p^{r}$-divisibility property concerns the map 
\(\namedright{\Omega\moore}{\Omega q}{\Omega S^{2n+1}}\). 
Since $\hlgy{\Omega S^{2n+1}}\cong\zmodp[x_{2n}]$, there 
is a natural filtration on $\Omega S^{2n+1}$ by skeleta. 
For $k\geq 0$, let \jk\ be the $2nk$-skeleton of $\Omega S^{2n+1}$. 
Then for each $k\geq 1$ there is a homotopy cofibration 
\[\lnameddright{S^{2nk-1}}{s_{k-1}}{\jkm}{}{\jk}.\] 
The filtation on $\Omega S^{2n+1}$ lets us put a filtration on  
$\Omega P^{2n+1}(p^{r})$. For $k\geq 0$ define the spaces $P_{k-1}$ 
and $B_{k-1}$, and the map $\gamma_{k-1}$, by the iterated homotopy 
pullback  
\[\diagram 
     B_{k-1}\rto\ddouble & P_{k-1}\rto^-{\gamma_{k-1}}\dto 
        & \Omega P^{2n+1}(p^{r})\dto^{\Omega q} \\ 
     B_{k-1}\rto & \jkm\rto & \Omega S^{2n+1}. 
  \enddiagram\] 
Since 
\(\namedright{\jkm}{}{\Omega S^{2n+1}}\) 
is a skeletal inclusion, a Serre spectral sequence calculation 
immediately implies that $B_{k-1}$ is $(2nk-2)$-connected and has 
a single cell in dimension $2nk-1$. Further, the composite 
\(S^{2nk-1}\hookrightarrow\namedright{B_{k-1}}{}{\jkm}\) 
is homotopic to $s_{k-1}$. Let $p_{k-1}$ be the composite 
\[p_{k-1}\colon S^{2nk-1}\hookrightarrow\namedright{B_{k-1}} 
     {}{P_{k-1}}.\] 
In Proposition~\ref{Pdivbyp} we show that the map $p_{k-1}$ is 
divisible by $p^{r}$ after suspending.  

\begin{proposition}
   \label{Pdivbyp}
   For $k\geq 1$, there is a homotopy commutative diagram
   \[\diagram
        \Sigma S^{2nk-1}\rto^-{\Sigma p_{k-1}}  
               \dto^-{\Sigma\underline{p}^{r}} 
           & \Sigma P_{k-1}\ddouble \\ 
        \Sigma S^{2nk-1}\rto^{\overline{p}_{k-1}} & \Sigma P_{k-1} 
     \enddiagram\]
   for some map $\overline{p}_{k-1}$.     
\end{proposition} 

\begin{proof}  
Since $\moore\simeq\Sigma P^{2n}(p^{r})$, the Bott-Samelson 
Theorem implies that 
$\hlgy{\Omega\moore}\cong T(\rhlgy{P^{2n}(p^{r})})\cong T(u,v)$, 
where $u$ and $v$ are in degrees $2n-1$ and $2n$ respectively, and 
the action of the Bockstein is determined by $\beta^{r}(v)=u$. 
In particular, $\beta^{r}(v^{k})\neq 0$. We will focus on
the pair $(v^{k},\beta^{r}(v^{k}))$.

Since $\moore\simeq\Sigma P^{2n}(p^{r})$, the James 
construction~\cite{J} implies that there is a homotopy equivalence
$\Sigma\Omega\moore\simeq
   \bigvee_{i=1}^{\infty}\Sigma (P^{2n}(p^{r}))^{(i)}$, 
where $(P^{2n})^{(i)}$ is the $i$-fold smash of $P^{2n}(p^{r})$ 
with itself. As the smash of two 
mod-$p^{r}$ Moore spaces is homotopy equivalent
to a wedge of two mod-$p^{r}$ Moore spaces, iterating shows that
$(P^{2n}(p^{r}))^{(i)}$ is homotopy equivalent to a wedge of mod-$p^{r}$ Moore
spaces for each $i\geq 2$, and hence so is $\Sigma\Omega\moore$. In 
particular, there is an inclusion of a wedge summand 
\(j\colon\namedright{P^{2nk+1}(p^{r})}{}{\Sigma\Omega\moore}\) 
whose image in homology is the pair $\{\sigma v^{k},\sigma\beta^{r}(v^{k})\}$. 
Define the space $R$ and the map $r$ by the homotopy cofibration 
\[\nameddright{P^{2np+1}(p^{r})}{j}{\Sigma\Omega\moore}{r}{R}.\] 

Define spaces $X$, $Y$, and $Z$ by the homotopy pullback diagram 
\[\diagram
     X\rto\ddouble & Y\rto\dto & Z\dto \\
     X\rto & \Sigma\pkm\rto^-{\Sigma\gamma_{k-1}}\dto^{r^{\prime}}
         & \Sigma\Omega\moore\dto^{r} \\
     & R\rdouble & R
  \enddiagram\]
where $r^{\prime}$ is defined as the composite 
$r\circ\Sigma\gamma_{k-1}$. We aim to identify $X$, $Y$, and $Z$ 
through a skeletal range. This will involve repeated applications 
of the Blakers-Massey Theorem, or equivalently, the Serre exact 
sequence. These state that if $A$ is $(m-1)$-connected and $C$ 
is $(n-1)$-connected then in dimensions~$\leq n+m-1$ the sequence 
\(\nameddright{A}{}{B}{}{C}\) 
is a homotopy fibration if and only if it is a homotopy cofibration. 

First, consider the homotopy fibration 
\(\nameddright{Z}{}{\Sigma\Omega\moore}{r}{R}\). 
By its definition, $R$ is $(2n-1)$-connected, so 
by the Blakers-Massey Theorem, the homotopy cofibration 
\(\nameddright{P^{2nk+1}(p^{r})}{j}{\Sigma\Omega\moore}{r}{R}\) 
is a homotopy fibration in dimensions~$\leq 2nk+2n-1$. But $Z$ 
is defined to be the homotopy fiber of $r$, so $Z\simeq P^{2nk+1}(p^{r})$ 
in dimensions~$\leq 2nk+2n-1$.  

Second, before considering $\Sigma\gamma_{k-1}$ 
we take another look at $\gamma_{k-1}$ and the homotopy fibration 
\(\nameddright{B_{k-1}}{}{P_{k-1}}{\gamma_{k-1}}{\Omega\moore}\). 
Recall that the $(2nk-1)$-skeleton of $B_{k-1}$ is $S^{2nk-1}$ and 
$p_{k-1}$ is the composite 
\(S^{2nk-1}\hookrightarrow\namedright{B_{k-1}}{}{P_{k-1}}\). 
Since $\Omega\moore$ is $(2n-2)$-connected, the Blakers-Massey Theorem 
implies that the sequence 
\(\nameddright{S^{2nk-1}}{p_{k-1}}{P_{k-1}}{\gamma_{k-1}} 
    {\Omega\moore}\) 
is a homotopy cofibration in dimensions~$\leq 2nk+2n-3$. Suspending, 
the sequence 
\(\llnameddright{S^{2nk}}{\Sigma p_{k-1}}{\Sigma P_{k-1}} 
    {\Sigma\gamma_{k-1}}{\Sigma\Omega\moore}\) 
is a homotopy cofibration in dimensions~$\leq 2nk+2n-1$. By the 
Blakers-Massey Theorem this sequence is also a homotopy fibration 
in dimensions~$\leq 2nk+2n-1$. On the other hand, the homotopy fiber 
of $\Sigma\gamma_{k-1}$ is defined to be $X$, so $X\simeq S^{2nk}$ in 
dimensions~$\leq 2nk+2n-1$, and the composite  
\(S^{2nk}\hookrightarrow\namedright{X}{}{\Sigma P_{k-1}}\)  
is homotopic to $\Sigma p_{k-1}$. 

Third, consider the homotopy fibration 
\(\nameddright{X}{}{Y}{}{Z}\). 
In dimensions~$\leq 2nk+2n-1$ we have $X\simeq S^{2nk}$ and 
$Z\simeq P^{2nk+1}(p^{r})$. The Blakers-Massey Theorem implies 
that this homotopy fibration is also a homotopy cofibration 
in this dimensional range. That is, in dimensions~$\leq 2nk+2n-1$, 
there is a homotopy cofibration 
\(\nameddright{S^{2nk}}{a}{Y}{}{P^{2nk+1}(p^{r})}\) 
for some map $a$. The corresponding long exact sequence in homology 
implies that $\hlgy{Y}\cong\hlgy{S^{2nk}}$ and $a_{\ast}=(p^{r})_{\ast}$ 
in dimension~$\leq 2nk+2n-1$. Thus $Y\simeq S^{2nk}$ and 
$a\simeq\underline{p}^{r}$ in dimensions~$\leq 2nk+2n-1$. 

Combining the conclusions of the previous two paragraphs gives a homotopy 
commutative diagram 
\[\diagram 
    S^{2nk}\rto^-{\underline{p}^{r}}\dto|<\hole|<<\ahook 
       & S^{2nk}\dto|<\hole|<<\ahook \\ 
    X\rto\ddouble & Y\dto \\ 
    X\rto & \Sigma P_{k-1} 
  \enddiagram\] 
in which the lower direction around the diagram is homotopic 
to $p_{k-1}$. The outer rectangle is therefore the diagram 
asserted by the Proposition.   
\end{proof} 

\begin{remark} 
When $k=1$ the divisibility property in Proposition~\ref{Pdivbyp} 
actually holds before suspending. This is seen directly from the 
homotopy fibration 
\(\nameddright{B_{0}}{}{P_{0}}{}{\Omega\moore}\) 
which, by the Blakers-Massey Theorem, is the homotopy cofibration 
\(\nameddright{S^{2n-1}}{\underline{p}^{r}}{S^{2n-1}}{} 
     {P^{2n}(p^{r})}\) 
in dimensions~$\leq 4n-3$.  
\end{remark}

\section{Extensions} 
\label{sec:extensions} 

In this section we prove Propositions~\ref{Eextension} 
and~\ref{Xextension}, which will complete the proof of 
Theorem~\ref{main}. 

\medskip\noindent 
\textit{Proof of Proposition~\ref{Eextension}}:  
Filter $E$ and $F$ by the spaces $E_{k}$ and $F_{k}$ considered in 
Section~\ref{sec:pdiv1}, so $E=\varinjlim E_{k}$ and $F=\varinjlim F_{k}$.  
By~(\ref{EkFk}), there are homotopy fibration sequences 
\(\namedddright{\Omega^{2} S^{2n+1}}{\partial_{k}}{E_{k}} 
     {f_{k}}{F_{k}}{}{\Omega S^{2n+1}}\). 
Note that when $k=0$ we have $F_{0}\simeq\ast$, so 
$E_{0}\simeq\Omega^{2} S^{2n+1}$ and $\partial_{0}$ is the 
identity map. Let  
\(\epsilon_{0}\colon\namedright{E_{0}}{}{\Omega^{2} S^{2np+1}\{p\}}\) 
be the map 
\(\namedright{\Omega^{2} S^{2n+1}}{S}{\Omega^{2} S^{2np+1}\{p\}}\).   
The asserted extension 
\(\namedright{E}{e_{1}}{\Omega^{2} S^{2np+1}\{p\}}\) 
of $S$ will be constructed as the limit of a sequence of extensions 
\(\namedright{E_{k}}{\epsilon_{k}}{\Omega^{2} S^{2np+1}\{p\}}\), 
where $\epsilon_{k}$ extends $\epsilon_{k-1}$.   
There are two cases: $k=1$ and $k>1$. 

The setup for both cases is the same. From the homotopy cofibration 
\(\nameddright{S^{2nk-1}}{g_{k-1}}{F_{k-1}}{}{F_{k}}\) 
and the map 
\(\namedright{E_{k}}{f_{k}}{F_{k}}\) 
we obtain a diagram of iterated homotopy pullbacks 
\[\diagram 
    \Omega^{2} S^{2n+1}\rdouble\dto 
      & \Omega^{2} S^{2n+1}\rdouble\dto^{\partial_{k}}  
       & \Omega^{2} S^{2n+1}\dto^{\partial_{k-1}} \\ 
    QE_{k-1}\rto\dto & E_{k-1}\rto\dto^{f_{k-1}} & E_{k}\dto^{f_{k}} \\ 
    S^{2nk-1}\rto^-{g_{k-1}} & F_{k-1}\rto & F_{k} 
  \enddiagram\] 
which defines the space $QE_{k-1}$. Lemma~\ref{DLpo} implies  
that $QE_{k-1}\simeq S^{2nk-1}\times\Omega^{2} S^{2n+1}$ and there 
is a homotopy pushout 
\begin{equation} 
  \label{Ekpushout} 
  \diagram 
    S^{2nk-1}\times\Omega^{2} S^{2n+1}\rto^-{\theta_{k-1}}\dto^{\pi_{2}} 
        & E_{k-1}\dto \\ 
    \Omega^{2} S^{2n+1}\rto^-{\partial_{k}} & E_{k}  
  \enddiagram 
\end{equation} 
for some map $\theta_{k-1}$. As the homotopy fibration 
\(\nameddright{\Omega^{2} S^{2n+1}}{\partial_{k-1}}{E_{k-1}} 
     {f_{k-1}}{F_{k-1}}\)  
is principal, by Remark~\ref{DLpoaction}, $\theta_{k-1}$ 
can be taken to be the composite 
\(\nameddright{S^{2nk-1}\times\Omega^{2} S^{2n+1}}{a\times 1} 
      {E_{k-1}\times\Omega^{2} S^{2n+1}}{\vartheta_{k-1}}{E_{k-1}}\), 
where $a$ is any choice of a lift of $g_{k-1}$ to $E_{k-1}$ 
and $\vartheta_{k-1}$ is the homotopy action from the principal 
fibration. 

Suppose $k=1$. Then $E_{0}\simeq\Omega^{2} S^{2n+1}$ and the 
action $\vartheta_{0}$ is the loop space multiplication. We 
claim that $a$ must be $E^{2}$, implying that 
$\theta_{0}\simeq E^{2}\cdot 1$. To see this, observe that in 
the homotopy fibration 
\(\nameddright{E_{1}}{f_{1}}{F_{1}}{}{\Omega S^{2n+1}}\) 
we have $F_{1}\simeq S^{2n}$ and the right map is degree one.  
Thus $E_{1}$ is $(4n-3)$-connected. 
Therefore the restriction of $\theta_{0}$ in~(\ref{Ekpushout}) 
must be degree one, that is, $a\simeq E^{2}$. Thus we have a 
homotopy pushout 
\[\diagram 
      S^{2n-1}\times\Omega^{2} S^{2n+1} 
           \rto^-{E^{2}\cdot 1}\dto^{\pi_{2}} 
         & \Omega^{2} S^{2n+1}\dto^{\partial_{1}} \\ 
      \Omega^{2} S^{2n+1}\rto^-{\partial_{1}} & E_{1}. 
  \enddiagram\]  
By Lemma~\ref{Ezerocase}, there is a homotopy commutative diagram 
\[\diagram 
      S^{2n-1}\times\Omega^{2} S^{2n+1} 
           \rto^-{E^{2}\cdot 1}\dto^{\pi_{2}} 
         & \Omega^{2} S^{2n+1}\dto^{S} \\  
      \Omega^{2} S^{2n+1}\rto^-{S} & \Omega^{2} S^{2np+1}\{p\}. 
  \enddiagram\] 
Hence there is a pushout map 
\(\epsilon_{1}:\namedright{E_{1}}{}{\Omega^{2} S^{2np+1}\{p\}}\) 
with the property that $\epsilon_{1}$ extends $\epsilon_{0}=S$. 

Next, suppose $k>1$ and assume that there is a map 
\(\epsilon_{k-1}:\namedright{E_{k-1}}{}{\Omega^{2} S^{2np+1}\{p\}}\) 
with the property that $\epsilon_{k-1}$ extends $\epsilon_{k-2}$. 
In~(\ref{Ekpushout}) we have seen that $\theta_{k-1}$ factors as 
\(\nameddright{S^{2nk-1}\times\Omega^{2} S^{2n+1}}{a\times 1} 
      {E_{k-1}\times\Omega^{2} S^{2n+1}}{\vartheta_{k-1}}{E_{k-1}}\), 
where $a$ lifts $g_{k-1}$. By Lemma~\ref{Edivbyp}, we can choose $a$ to 
be $\overline{g}_{k-1}$, which has the property that  
$\overline{g}_{k-1}\simeq\overline{h}_{k-1}\circ\underline{p}^{r}$ 
for some map $\overline{h}_{k-1}$. This fulfils one hypothesis of 
Theorem~\ref{extension}. The other is fulfilled because the $H$-space 
$\Omega^{2} S^{2np+1}\{p\}$ is homotopy associative and by~\cite{N3} 
its $p^{th}$-power map is null homotopic. Thus Theorem~\ref{extension} 
implies that there is an extension  
\[\diagram 
      E_{k-1}\rto^-{\epsilon_{k-1}}\dto 
        & \Omega^{2} S^{2np+1}\{p\}\ddouble \\ 
      E_{k}\rto^-{\epsilon_{k}} & \Omega^{2} S^{2np+1}\{p\}   
  \enddiagram\] 
for some map $\epsilon_{k}$. 
$\qqed$  

As described in~(\ref{l1dgrm}), the extension in 
Proposition~\ref{Eextension} implies there is a homotopy fibration sequence 
\(\namedddright{\Omega^{2} S^{2n+1}}{}{Y}{}{X}{q_{X}}{\Omega S^{2n+1}}\),  
and in Proposition~\ref{Xextension} we want to extend the map 
\(\namedright{Y}{g}{W_{np}}\) 
to a map 
\(\namedright{X}{}{W_{np}}\). 
To do this we will filter $X$ and produce iterated extensions using 
Theorem~\ref{extension}, but first a certain decomposition issue needs 
to be resolved. 
 
For $k\geq 0$, define the space $X_{k}$ and the map $q_{k}$ by the 
homotopy pullback diagram 
\[\diagram 
       X_{k}\rto\dto^{q_{k}} & X\dto^{q_{X}} \\ 
       \jk\rto & \Omega S^{2n+1}. 
  \enddiagram\] 
From the homotopy cofibration 
\(\nameddright{S^{2nk-1}}{s_{k-1}}{\jkm}{}{\jk}\) 
and the map 
\(\namedright{X_{k}}{}{\jk}\),  
we obtain a diagram of iterated homotopy pullbacks  
\begin{equation} 
   \label{Xkdgrm} 
   \diagram 
       Y\rdouble\dto & Y\rdouble\dto & Y\dto \\ 
       QX_{k-1}\rto\dto & X_{k-1}\rto\dto^{q_{k-1}} 
          & X_{k}\dto^{q_{k}} \\ 
       S^{2nk-1}\rto^-{s_{k-1}} & \jkm\rto & \jk 
   \enddiagram 
\end{equation} 
which defines the space $QX_{k-1}$. Lemma~\ref{DLpo} implies that 
$QX_{k-1}\simeq S^{2nk-1}\times Y$ and there is a homotopy pushout 
\[\diagram 
    QX_{k-1}\simeq S^{2nk-1}\times Y\rto^-{\theta_{k-1}}\dto^{\pi_{2}} 
        & X_{k-1}\dto \\ 
    Y\rto^-{\delta_{k}} & X_{k}  
  \enddiagram\] 
for some map $\theta_{k-1}$. However, the fibration 
\(\nameddright{Y}{}{X_{k-1}}{q_{k-1}}{\jkm}\) 
may not be principal so we cannot appeal to Remark~\ref{DLpoaction} 
in order to show that the equivalence for $QX_{k-1}$ can be chosen 
so that $\theta_{k-1}$ factors in the manner required by 
Theorem~\ref{extension}. To obtain such a factorization we have to 
work a bit harder, and this is the purpose of the next three lemmas. 

Recall from~(\ref{OmegaPXdgrm}) that the  map  
\(\namedright{\Omega P^{2n+1}(p^{r})}{\Omega q}{\Omega S^{2n+1}}\) 
factors as the composite 
\(\namedddright{\Omega P^{2n+1}(p^{r})}{i_{X}}{X}{} 
    {\Omega S^{2n+1}\{p^{r}\}}{}{\Omega S^{2n+1}}\). 
So for each $k\geq 1$ there is an iterated homotopy pullback 
diagram  
\begin{equation} 
  \label{XkCk} 
  \diagram 
     P_{k-1}\rto\dto & X_{k-1}\rto\dto & C_{k-1}\rto\dto 
        & \jkm\dto \\ 
     \Omega P^{2n+1}(p^{r})\rto^-{i_{X}} & X\rto  
        & \Omega S^{2n+1}\{p^{r}\}\rto & \Omega S^{2n+1}    
  \enddiagram 
\end{equation}  
which defines the space $C_{k-1}$. Observe that the left and middle 
homotopy pullbacks in~(\ref{XkCk}) give a rectangle which is a homotopy 
pullback. In particular, this implies that there is a homotopy fibration 
\(\nameddright{P_{k-1}}{}{C_{k-1}}{}{E}\). 
Also, the right homotopy pullback implies that there is a homotopy 
fibration 
\(\nameddright{C_{k-1}}{}{\jkm}{j_{k-1}}{\Omega S^{2n+1}}\), 
where $j_{k-1}$ is the composite 
\(j_{k-1}\colon\jkm\hookrightarrow\namedright{\Omega S^{2n+1}}{p^{r}} 
     {\Omega S^{2n+1}}\).  
In particular, there is a homotopy action 
\(\namedright{C_{k-1}\times\Omega^{2} S^{2n+1}}{}{C_{k-1}}\). 
Let $u$ be the composite 
\[u:\nameddright{P_{k-1}\times\Omega^{2} S^{2n+1}}{} 
     {C_{k-1}\times\Omega^{2} S^{2n+1}}{}{C_{k-1}}.\]  

\begin{lemma} 
   \label{Yaction1} 
   There is a map 
   \(\namedright{P_{k-1}\times Y}{}{X_{k-1}}\) 
   and a homotopy pullback 
   \[\diagram 
        P_{k-1}\times Y\rto\dto & X_{k-1}\dto \\ 
        P_{k-1}\times\Omega^{2} S^{2n+1}\rto^-{u} & C_{k-1}. 
     \enddiagram\] 
\end{lemma} 

\begin{proof} 
Consider the diagram of extended homotopy pullbacks  
\[\diagram 
     \Omega^{2} S^{2n+1}\rto\ddouble & C_{k-1}\rto\dto 
        & \jkm\rto^-{j_{k-1}}\dto & \Omega S^{2n+1}\ddouble \\ 
     \Omega^{2} S^{2n+1}\rto\ddouble & \Omega S^{2n+1}\{p^{r}\}\rto\dto 
        & \Omega S^{2n+1}\rto^-{p^{r}}\dto & \Omega S^{2n+1}\ddouble \\ 
     \Omega^{2} S^{2n+1}\rto & E\rto & F\rto & \Omega S^{2n+1}. 
  \enddiagram\] 
Here, the lower ladder is from the definitions of the spaces $E$ and $F$ 
in~(\ref{EFdgrm}) and the upper ladder is from the definition of the 
space $C_{k-1}$. The naturality of the canonical homotopy action for 
a fibration sequence implies that there is a homotopy commutative diagram 
of actions 
\begin{equation} 
  \label{actiondgrm} 
  \diagram 
     C_{k-1}\times\Omega^{2} S^{2n+1}\rto\dto & C_{k-1}\dto \\ 
     \Omega S^{2n+1}\{p^{r}\}\times\Omega^{2} S^{2n+1}\rto\dto 
         & \Omega S^{2n+1}\{p^{r}\}\dto \\ 
     E\times\Omega^{2} S^{2n+1}\rto & E. 
  \enddiagram 
\end{equation}  
Combining this with the homotopy fibration 
\(\nameddright{P_{k-1}}{}{C_{k-1}}{}{E}\) 
gives a homotopy commutative diagram 
\begin{equation} 
  \label{PkCkaction} 
  \diagram 
      P_{k-1}\times\Omega^{2} S^{2n+1}\rto\dto 
          & C_{k-1}\times\Omega^{2} S^{2n+1}\rto\dto & C_{k-1}\dto \\ 
      \ast\times\Omega^{2} S^{2n+1}\rto 
          & E\times\Omega^{2} S^{2n+1}\rto & E. 
  \enddiagram 
\end{equation}  
Note that the upper row of this diagram is the definition of $u$. 
Now compose the diagram with the map 
\(\namedright{E}{e_{1}}{\Omega^{2} S^{2np+1}\{p\}}\). 
Let $v$ and $w$ be the composites 
\[v:\nameddright{C_{k-1}}{}{E}{e_{1}}{\Omega^{2} S^{2np+1}\{p\}}\]  
\[w:\nameddright{P_{k-1}\times\Omega^{2} S^{2n+1}}{u}{C_{k-1}}{v} 
     {\Omega^{2} S^{2np+1}\{p\}}.\]  
We wish to identify the homotopy fibers of $v$ and $w$. First, $v$ 
factors as the composite 
\(\namedddright{C_{k-1}}{}{\Omega S^{2n+1}\{p^{r}\}}{}{E}{e_{1}} 
    {\Omega^{2} S^{2np+1}\{p\}}\). 
The latter two maps in this composite define $\epsilon_{1}$, whose 
homotopy fiber is $X$. Thus the middle pullback in~(\ref{XkCk}) 
implies that the homotopy fiber of $v$ is~$X_{k-1}$. Second, the 
homotopy commutativity of~(\ref{PkCkaction}) implies that $w$ is 
homotopic to the composite 
\(\namedddright{P_{k-1}\times\Omega^{2} S^{2n+1}}{\pi_{2}} 
    {\Omega^{2} S^{2n+1}}{}{E}{e_{1}}{\Omega^{2} S^{2np+1}\{p\}}\), 
which in turn is homotopic to the composite 
\(\nameddright{P_{k-1}\times\Omega^{2} S^{2n+1}}{\pi_{2}} 
    {\Omega^{2} S^{2n+1}}{S}{\Omega^{2} S^{2np+1}\{p\}}\). 
As the homotopy fiber of $S$ is $Y$, this implies that there 
is a homotopy fibration 
\(\nameddright{P_{k-1}\times Y}{}{P_{k-1}\times\Omega^{2} S^{2n+1}} 
     {w}{\Omega^{2} S^{2np+1}\{p\}}\). 
Putting the fibrations for $v$ and $w$ together, we obtain a homotopy 
pullback diagram 
\[\diagram 
     P_{k-1}\times Y\rto\dto & X_{k-1}\dto \\ 
     P_{k-1}\times\Omega^{2} S^{2n+1}\rto^-{u}\dto^{w} 
        & C_{k-1}\dto^{v} \\ 
     \Omega^{2} S^{2np+1}\{p\}\rdouble & \Omega^{2} S^{2np+1}\{p\}  
  \enddiagram\] 
which proves the lemma. 
\end{proof} 

Define $\overline{\theta}_{k-1}$ as the composite 
\[\overline{\theta}_{k-1}\colon\llnameddright{S^{2nk-1}\times Y} 
    {p_{k-1}\times 1}{P_{k-1}\times Y}{}{X_{k-1}}.\] 

\begin{lemma} 
   \label{Yaction2} 
   There is a homotopy pullback 
   \[\diagram 
        S^{2nk-1}\times Y\rto^-{\overline{\theta}_{k-1}}\dto^{\pi_{1}} 
           & X_{k-1}\dto^{q_{k-1}} \\ 
        S^{2nk-1}\rto^-{s_{k-1}} & \jkm. 
     \enddiagram\] 
\end{lemma} 

\begin{proof} 
Consider the diagram 
\begin{equation} 
  \label{QXtrivdgrm} 
  \diagram 
     S^{2nk-1}\times Y\rto^-{p_{k-1}\times 1}\ddto^{\pi_{1}} 
        & P_{k-1}\times Y\rto\dto & X_{k-1}\dto \\  
     & P_{k-1}\times\Omega^{2} S^{2n+1}\rto^-{u}\dto^{\pi_{1}} 
        & C_{k-1}\dto \\ 
     S^{2nk-1}\rto^-{p_{k-1}} & P_{k-1}\rto & \jkm. 
  \enddiagram 
\end{equation}  
The left rectangle is a homotopy pullback by the naturality of the 
projection. The upper right square is a homotopy pullback by 
Lemma~\ref{Yaction1}. For the lower right square, by definition, 
$u$~is the composite  
\(\nameddright{P_{k-1}\times\Omega^{2} S^{2n+1}}{} 
    {C_{k-1}\times\Omega^{2} S^{2n+1}}{}{C_{k-1}}\) 
where the right map is the action associated to the homotopy fibration 
\(\nameddright{C_{k-1}}{}{\jkm}{j_{k-1}}{\Omega S^{2n+1}}\). 
A canonical property of such an action is that there is a  
homotopy commutative diagram 
\[\diagram 
      C_{k-1}\times\Omega^{2} S^{2n+1}\rto\dto^{\pi_{1}} 
         & C_{k-1}\dto \\ 
      C_{k-1}\rto & \jkm. 
  \enddiagram\] 
Note that this is in fact a homotopy pullback. When precomposed  
with the map 
\(\namedright{P_{k-1}\times\Omega^{2} S^{2n+1}}{} 
    {C_{k-1}\times\Omega^{2} S^{2n+1}}\), 
the naturality of the projection $\pi_{1}$ implies that we obtain 
the lower right square in~(\ref{QXtrivdgrm}), and that this square is 
a homotopy pullback. Thus~(\ref{QXtrivdgrm}) consists of homotopy 
pullbacks and so its outer perimeter is also a homotopy pullback. 
As the top row of~(\ref{QXtrivdgrm}) is the definition 
of $\overline{\theta}_{k-1}$ and its right column is $q_{k-1}$, 
the outer perimeter of~(\ref{QXtrivdgrm}) is the diagram 
asserted by the lemma.  
\end{proof} 

Lemma~\ref{Yaction2} describes the homotopy pullback of 
\(\namedright{X_{k-1}}{}{\jkm}\) 
and 
\(\namedright{S^{2nk-1}}{s_{k-1}}{\jkm}\). 
But by~(\ref{Xkdgrm}), $QX_{k-1}$ is defined to be this pullback. 
Thus there is a homotopy equivalence 
\(\namedright{S^{2nk-1}\times Y}{\simeq}{QX_{k-1}}\) 
with an appropriate action property, as stated in the following. 

\begin{lemma} 
   \label{Xdivbyp} 
   There is a decomposition 
   $QX_{k-1}\simeq S^{2nk-1}\times Y$ and a homotopy pullback 
   \[\diagram 
        QX_{k-1}\simeq S^{2nk-1}\times Y
              \rto^-{p_{k-1}\times 1}\dto^{\pi_{1}} 
           & P_{k-1}\times Y\rto & X_{k-1}\dto \\ 
        S^{2nk-1}\rrto^-{s_{k-1}} & & \jkm   
     \enddiagram\]  
   where $p_{k-1}$ is the map in Proposition~\ref{Pdivbyp}. 
   In particular, for $k\geq 1$ we have 
   $\Sigma p_{k-1}\simeq\overline{p}_{k-1}\circ\underline{p}^{r}$ 
   for some map $\overline{p}_{k-1}$.~$\qqed$ 
\end{lemma} 

We are now ready to prove the existence of the second extension. 

\medskip\noindent 
\textit{Proof of Proposition~\ref{Xextension}}:  
Filter $X$ and $\Omega S^{2n+1}$ by the spaces $X_{k}$ and \jk,  
so $X=\varinjlim X_{k}$ and $\Omega S^{2n+1}=\varinjlim\jk$.  
By~(\ref{Xkdgrm}), there are homotopy fibrations   
\(\nameddright{Y}{\delta_{k}}{X_{k}}{q_{k}}{\jk}\). 
Note that when $k=0$ we have $J_{0}(S^{2n})\simeq\ast$, so 
$X_{0}\simeq Y$ and $\delta_{0}$ is the identity map. Let  
\(\varepsilon_{0}\colon\namedright{X_{0}}{}{W_{np}}\) 
be the map 
\(\namedright{Y}{g}{W_{np}}\).  
The asserted extension 
\(\namedright{X}{e_{2}}{W_{np}}\) 
of $g$ will be constructed as the limit of a sequence of extensions 
\(\namedright{X_{k}}{\varepsilon_{k}}{\Omega W_{np}}\) 
with the property that $\varepsilon_{k}$ extends $\varepsilon_{k-1}$.  

If $k=1$ we have 
\(\namedright{X_{0}}{\varepsilon_{0}}{\Omega W_{np}}\) 
in place. If $k>1$ assume that there is a map 
\(\varepsilon_{k-1}\colon\namedright{X_{k-1}}{} 
     {\Omega W_{np}}\) 
with the property that $\varepsilon_{k-1}$ is an extension 
of $\varepsilon_{k-2}$. We set up to use Theorem~\ref{extension}. 
From the homotopy cofibration 
\(\nameddright{S^{2nk-1}}{s_{k-1}}{\jkm}{}{\jk}\) 
and the map 
\(\namedright{X_{k}}{q_{k}}{\jk}\) 
we obtain a diagram of iterated homotopy pullbacks  
\[\diagram 
     Y\rdouble\dto & Y\rdouble\dto^{\delta_{k-1}} 
        & Y\dto^{\delta_{k}} \\ 
     QX_{k-1}\rto\dto & X_{k-1}\rto\dto^{q_{k}} 
        & X_{k}\dto^{q_{k-1}} \\ 
     S^{2nk-1}\rto^-{s_{k-1}} & \jkm\rto & \jk  
  \enddiagram\]  
By Lemma~\ref{DLpo} there is a decomposition  
$QX_{k-1}\simeq S^{2nk-1}\times Y$ and a homotopy pushout 
\[\diagram 
    QX_{k-1}\simeq S^{2nk-1}\times Y\rto^-{\theta_{k-1}}\dto^{\pi_{2}} 
        & X_{k-1}\dto \\ 
    Y\rto^-{\delta_{k}} & X_{k}  
  \enddiagram\] 
for some map $\theta_{k-1}$. By hypothesis, there is a map 
\(\namedright{X_{k-1}}{\varepsilon_{k-1}}{W_{np}}\). 
Now we check that the two hypotheses of Theorem~\ref{extension} 
are satisfied. First, by Lemma~\ref{Xdivbyp}, the map  
\(\namedright{S^{2nk-1}\times Y}{\theta_{k-1}}{X_{k-1}}\) 
in this pushout can be taken to be the composite 
\(\llnameddright{S^{2nk-1}\times Y}{p_{k-1}\times 1}
    {P_{k-1}\times Y}{}{X_{k-1}}\),  
where $\Sigma p_{k-1}\simeq\overline{t}_{k-1}\circ\underline{p}^{r}$ 
for some map $\overline{p}_{k-1}$. Second, as $W_{np}$ is a loop 
space it is homotopy associative, and by Lemma~\ref{BWnpexp} its 
$p^{th}$-power map is null homotopic. Thus Theorem~\ref{extension} 
implies that there is an extension  
\[\diagram 
      X_{k-1}\rto^-{\varepsilon_{k-1}}\dto 
        & W_{np}\ddouble \\ 
      X_{k}\rto^-{\varepsilon_{k}} & W_{np}   
  \enddiagram\] 
for some map $\varepsilon_{k}$. 
$\qqed$

\section{Consequences} 
\label{sec:consequences} 

In this section we explore some consequences of our proof of 
Theorem~\ref{main} in order to round off the picture. First, 
in Lemma~\ref{BWnpexp} we used a factorization of the 
$p^{th}$-power map on $\Omega^{2} S^{2np+1}$ through the double 
suspension to show that the $p^{th}$-power map on $BW_{np}$ is 
null homotopic. Now that we know the $p^{th}$-power map on 
$\Omega^{2} S^{2n+1}$ factors through the double suspension we 
can use the same argument to improve from the $np$-case to all 
cases. 

\begin{theorem} 
   \label{Bwnexp} 
   The $p{th}$-power map on $BW_{n}$ is null homotopic.~$\qqed$ 
\end{theorem} 

Now return to the construction of the extension 
\[\diagram 
     \Omega^{2} S^{2n+1}\rto^-{\partial_{E}}\dto^{S} 
         & E\dlto^-{e_{1}} \\ 
     \Omega^{2} S^{2np+1}\{p\} & 
  \enddiagram\] 
in Proposition~\ref{Eextension}. The only properties of 
$\Omega^{2} S^{2np+1}\{p\}$ and $S$ used in the proof were: 
(i) $S$ is an $H$-map, (ii) $\Omega^{2} S^{2np+1}\{p\}$ is a 
homotopy associative $H$-space, and (iii) the $p^{th}$-power 
map on $\Omega^{2} S^{2np+1}\{p\}$ is null homotopic. Recall 
that the map 
\(\namedright{\Omega^{2} S^{2n+1}}{\nu}{BW_{n}}\) 
is an $H$-map. By~\cite{G} for $p\geq 5$ and~\cite{T2} for $p=3$, 
$BW_{n}$ is a homotopy associative $H$-space. Thus, now that we 
know the $p^{th}$-power map on $BW_{n}$ is null homotopic we can 
replace $\Omega^{2} S^{2np+1}\{p\}$ by $BW_{n}$ in the proof of 
Proposition~\ref{Eextension} to obtain the following. 

\begin{proposition} 
   \label{BWnextension} 
   There is an extension 
   \[\diagram 
         \Omega^{2} S^{2n+1}\rto^-{\partial_{E}}\dto^{\nu} 
             & E\dlto^-{e} \\ 
         BW_{n} & 
     \enddiagram\] 
   for some map $e$.~$\qqed$ 
\end{proposition} 

Proposition~\ref{BWnextension} recovers the key result 
in~\cite{GT}. It allows for a reformulation of \anick\ which 
satisfies additional properties. Since $\partial_{E}$ factors 
as the composite 
\(\nameddright{\Omega^{2} S^{2n+1}}{}{\Omega\curly}{}{E}\), 
there is a homotopy pullback diagram  
\[\diagram 
    S^{2n-1}\rto\dto^{E^{2}} & \anick\rto\dto 
        & \Omega S^{2n+1}\ddouble \\  
    \Omega^{2} S^{2n+1}\rto\dto^{\nu} & \Omega\curly\rto\dto 
        & \Omega S^{2n+1} \\ 
    BW_{n}\rdouble & BW_{n} &   
  \enddiagram\] 
which reformulates the construction of \anick. Further, there 
is a homotopy pullback diagram 
\[\diagram   
      \Omega P^{2n+1}(p^{r})\rto\ddouble 
        & \anick\rto\dto & X\dto \\  
      \Omega P^{2n+1}(p^{r})\rto^-{\Omega i}  
        & \Omega\curly\rto\dto & E\dto \\  
      & BW_{n}\rdouble & BW_{n}  
  \enddiagram\] 
which implies the additional property that $\Omega i$ factors 
through \anick. Carrying on, a key property of~\cite{CMN1} regarding 
the homotopy theory of the Moore space can be recovered, namely, 
that the inclusion 
\(\imath\colon\namedright{S^{2n-1}}{}{\Omega F}\) 
of the bottom cell has a left homotopy inverse. Since $\Omega i$ 
factors through \anick, so does $\Omega q$. Thus there is a 
homotopy pullback diagram  
\[\diagram  
     \Omega F\rto\dto^{r} & \Omega P^{2n+1}(p^{r})\rto^-{\Omega q}\dto 
         & \Omega S^{2n+1}\ddouble \\ 
     S^{2n-1}\rto & \anick\rto & \Omega S^{2n+1}  
  \enddiagram\] 
which defines the map $r$. By connectivity, the entire left square is 
degree one in $H_{2n-1}(\ )$, and so $r\circ\imath$ is degree one, 
implying that it is homotopic to the identity map.

%%% The bibliography %%% 
\bibliographystyle{amsalpha}

\end{document}